\documentclass[12pt]{article}

      \usepackage{latexsym}
         \usepackage[reqno, namelimits, sumlimits]{amsmath}
         \usepackage{amssymb, amsfonts}
         \usepackage{amsthm}

 \newtheorem{theorem}{Theorem}[section]
 
 \newtheorem{definition}[theorem]{Definition}
 
 \newtheorem{lemma}[theorem]{Lemma}
 
 \newtheorem{remark}[theorem]{Remark}

 \newtheorem{pro}[theorem]{Proposition}
 \newtheorem{con}[theorem]{Conjecture}
\title{Rescalings at possible singularities of Navier-Stokes equations in half space
}

 \author{G. Seregin,  V. \v Sver\'ak 
 }

\begin{document}
\maketitle

\setcounter{equation}{0}
\section{Introduction}

We consider the  initial boundary-value problem for the incompressible Navier-Stokes equations in half-space $\mathbb R^3_+=\{x_3>0\}$:
\begin{eqnarray}
\left.\begin{array}{rcl}
v_t+v\nabla v -\Delta v & = & -\nabla q \\
{\rm div \,}v & = & \,\,\,\,\,0
\end{array} \right\}& &  \quad \hbox{in $\mathbb R^3_+\times(0,\infty)$}\,, \label{mainequation} \\
v(\,\cdot\,,t)|_{\partial \mathbb R^3_+}\,  =   \,\,\,\,\,\,\,\,0\,\,\,\,\,\,  \quad  & & \quad \hbox{for $t>0$}\,,\label{mbcv} \\
v(\,\cdot\,,0) \,\, = \,\,\,\,\,\,\,v_0\,\,\, \,\,\quad & & \quad \hbox{in $\mathbb R^3$}\,,\label{micv}
\end{eqnarray}
where $v_0$ is a sufficiently regular div-free field in $\mathbb R^3_+$ with sufficiently fast decay  for $x\to\infty$ and $v_0|_{\partial \mathbb R^3_+}=0\,$.

Our main goal is to understand connections between  possible blow-up of strong solutions and Liouville theorem for bounded ancient mild solutions, in the spirit considered for all space in~\cite{KochNadSS2009}.  We recall that the local-in-time existence of strong solutions was proved in all space in \cite{Le}  and for bounded domains similar results appeared in \cite{LadKis}. As to  unbounded domains, we refer, for example, to \cite {L} or \cite{Heywood1978}.

Let us consider the local-in-time strong solution $v$ on its maximal interval of existence $[0,T)$, which will be assumed to be finite.\footnote{At the time of this writing it is unknown whether such solutions exist.} The time $T$ is then the blow-up time. Such solutions are known to be unique (unlike the weak solutions, which on the other hand are global). We set
\begin{equation}\label{blowup}
    g(t)=\sup\limits_{0<\tau\leq t}M(\tau)\,\,,
\end{equation}
where
$$M(t)=\sup\limits_{x\in \mathbb R^3_+}|v(x,t)|\,.$$
It is well-known that $g(t)\to+\infty$ as $t\to T_-$, see for example~\cite{Le}.

Our considerations are motivated by methods used in the theory of geometric flows and based on scale-invariant transformations of the solution $v$ when time is approaching $T$, see~\cite{KochNadSS2009} for a more detailed discussion and references.  In the case of the Navier-Stokes equations, the scale-invariant transformations have the form $v(x,t), q(x,t)\to \lambda v(\lambda x, \lambda^2t)\,,\, \lambda^2 q(\lambda x, \lambda^2t)$. In the whole space,  taking  limits of sequences of suitably scaled solutions produces  so-called bounded ancient (backward) solutions to the Navier-Stokes equations.  They are defined on the semi-infinite time interval $]-\infty,0]$ (backward in time), are bounded in $L^\infty$  and  in fact belong to a subclass called mild bounded ancient solutions. In \cite{KochNadSS2009}, it has been shown that mild bounded ancient solutions are infinitely smooth both in space and in time. Moreover, given that $T$ is a blowup time, the limiting mild bounded ancient solution cannot vanish. It has been conjectured in \cite{KochNadSS2009} that \textit{any mild bounded ancient solution is a constant}. This would rule out blow-ups of Type I in the case of the Cauchy problem (whole space). Let us recall that Type I blowup is usually defined by the inequality
\begin{equation}\label{type1wholespace}
    |v(x,t)|\leq \frac c{\sqrt{T-t}}
\end{equation}
for any $x\in\mathbb R^3$ and any $t<T$, although more general definitions are also possible.

To show a difference between bounded ancient solutions in the whole space and in the half space, we recall their definitions. We say that $u$ is a \textit{bounded ancient solution} of the Navier-Stokes equations if $u$ is  bounded in $Q_-=\mathbb R^3\times ]-\infty,0[$
and satisfies these equations in the sense of distributions with divergence free test functions, i.e.,
\begin{equation}\label{momentum}
    \int\limits_{Q_-}\Big(u\cdot(\partial_t\varphi+\Delta \varphi)+u\otimes u:\nabla\varphi\Big)dz=0
\end{equation}
for any $\varphi\in C^\infty_{0,0}(Q_-):=\{\varphi\in C^\infty_{0}(Q_-):\,\,\mbox{div}\,\varphi=0\}$;
\begin{equation}\label{incompressibility}
    \int\limits_{Q_-}u\cdot\nabla q dz=0
\end{equation}
for any $q\in C^\infty_0(Q_-)$. By scaling, we may assume that $|u|$ is bounded by one.
There are simple nontrivial bounded ancient solutions of the form
\begin{equation}\label{nontrivialsolutioninwhole}
 u(x,t)=a(t)
\end{equation}
where $a$ is an arbitrary bounded function of $t$ only.

A vector field $u$ is called a \textit{mild bounded ancient solution} if $u$ is a bounded ancient solution and there exists a pressure field $p\in L_\infty(-\infty,0;BMO(\mathbb R^3))$ such that
\begin{equation}\label{momentumwithpressure}
    \int\limits_{Q_-}\Big(u\cdot(\partial_t\varphi+\Delta \varphi)+u\otimes u:\nabla\varphi\Big)dz=-\int\limits_{Q_-}p\,{\rm div}\,\varphi dz
\end{equation}
for any $\varphi\in C^\infty_{0}(Q_-)$. It is not so difficult to see that any solution of the form (\ref{nontrivialsolutioninwhole}) is a mild bounded ancient solution if and only if $a(t)=constant$. As already mentioned above, mild bounded ancient solutions are infinitely smooth and have the following property: for any $A<0$, they can be presented in the form
$$u_i(x,t)=\int\limits_{\mathbb R^3}\Gamma(x-y,t-A)u_i(y,A)dy+$$
\begin{equation}\label{mildinwholespace}
+\int\limits^t_A \int\limits_{\mathbb R^3} K_{ijm}(x,y,t-\tau)u_j(y,\tau)u_m(y,\tau)dy d\tau,
\end{equation}
where $\Gamma$ is the well known heat kernel and $K$ is obtained from the Oseen tensor by differentiation in spatial variables, see details in \cite{Sol2003} and \cite{KochNadSS2009}.  That is why those solutions are called mild ones. By the way, it is the original definition of mild bounded ancient solutions given in \cite{KochNadSS2009} and the equivalent definition in terms of the pressure appeared in \cite{SerSve2009} later on.

The case of the half space is more complicated. First, it is not immediately clear how to understand the homogeneous Dirichlet boundary conditions when the velocity $u$ is only bounded in $Q_-^+:=\{z=(x,t):\,\,x\in\mathbb R^3_+,\,\,-\infty<t<0\}$. We shall use the following weak definition of~(\ref{mbcv}):
\begin{equation}\label{momentum in hs}
    \int\limits_{Q_-^+}\Big(u\cdot(\partial_t\varphi+\Delta \varphi)+u\otimes u:\nabla\varphi\Big)dz=0
\end{equation}
for any $\varphi\in C^\infty_{0,0}(Q_-)$ with $\varphi(x',0,t)=0$ for any $x'\in\mathbb R^2$ and for any $-\infty<t<0$;
\begin{equation}\label{incompressibility in hs}
    \int\limits_{Q_-^+}u\cdot\nabla q dz=0
\end{equation}
for any $q\in C^\infty_0(Q_-)$.

In \cite{JiaSS2012}, a class of simple non-trivial bounded ancient solutions to the Navier-Stokes equations has been presented. They describe a bounded shear flow in the half space and have the form
\begin{equation}\label{boundedshearflow}
    u(x,t)=(u_1(x_3,t),u_2(x_3,t),0).
\end{equation}
It has also been shown in \cite{JiaSS2012} that there are no other non-trivial solutions to the linear Stokes system in the half space. It is unknown whether or not this is true for the nonlinear case.

We now  define mild bounded ancient solutions in a half space.
\begin{definition}\label{hsmildsol}  A bounded function $u$ is a mild bounded ancient solution if and only if  there exists
a pressure $p$ such that $p=p^1+p^2$,
where the even extension of $p^1$ to the whole $\mathbb R^3$ with respect to $x_3$ is $L_\infty(-\infty,0;BMO(\mathbb R^3))$-function,
\begin{equation}\label{p1}
 \triangle p^1=-{\rm div div}\,u\otimes u
\end{equation}
in $Q^+_-$ with $p^1_{,3}(x',0,t)=0$  and $p^2(\cdot,t)$ is a harmonic function in $\mathbb R^3_+$ whose gradient satisfies the estimate
\begin{equation}\label{log}
  |\nabla p^2(x,t)|\leq c\ln(2+ 1/{x_3})
\end{equation}
for all $(x,t)\in Q_-^+$
and has the property
 \begin{equation}\label{p2}
   \sup\limits_{x'\in\mathbb R^2}|\nabla p^2(x,t)|\to 0
 \end{equation}
as $x_3\to \infty$;  $u$ and $p$ satisfy (\ref{incompressibility in hs}) and
\begin{equation}\label{momentumpres}
    \int\limits_{Q^+_-}\Big(u\cdot(\partial_t\varphi+\Delta \varphi)+u\otimes u:\nabla\varphi+p{\rm div}\,\varphi\Big)dx dt=0
\end{equation}
for any $\varphi\in C^\infty_{0}(Q_-)$ with $\varphi(x',0,t)=0$ for any $x'\in \mathbb R^2$ and for any $t<0$.
\end{definition}
 \begin{remark}\label{smoothms} If $u$ is a mild bounded ancient solution, then $\nabla u\in L_\infty(Q_-^+)$. The function $u$ is infinitely smooth in spatial variables in upper half space $x_3>0$. \end{remark}

As  in the case of the whole space, mild bounded ancient solutions  can be defined by formula (\ref{mildinwholespace}),  in which $\mathbb R^3$ is replaced with $\mathbb R^3_+$ and kernels with their half-space analogues, for details see Section 2. The corresponding statement might be called the equivalence theorem. The proof of such a result is more involved than
its whole space version and will be published elsewhere.

It is worth noticing that a nontrivial solution of the form (\ref{boundedshearflow}) is a mild bounded ancient solutions if and only $u=0$. Indeed, 
for solutions
(\ref{boundedshearflow}), the gradient of the pressure is a function of time only. By the above theorem, this is only possible if the gradient of the pressure is equal to zero. In turn, this means that each component $u_\alpha$, $\alpha=1,2$ is a bounded ancient solution to the heat equation in $\{x_3>0\}\times ]-\infty,0[$ with the boundary condition $u_\alpha(0,t)=0$, which implies that $u_\alpha\equiv0$.


We believe that the following is true:
\begin{con}\label{conjecture} There is no non-trivial mild bounded ancient solution to the Navier-Stokes equations in the half space.\end{con}

The validity of Conjecture \ref{conjecture} and the conjecture, made in \cite{KochNadSS2009} and mentioned above, would rule out Type I blowups in the broad sense, when understood as blow-up solutions with a suitable bounded scale-invariant quantity.

Let us state our main result.
\begin{theorem}\label{main}
Assume that initial boundary value problem (\ref{mainequation})-(\ref{micv}) has a solution that blows up at time $T$. There exists at least one non-trivial (non-zero) mild bounded ancient solution either in the whole space  or in the half space. \end{theorem}

 The appearance of mild bounded ancient solutions in the whole space is not surprising and it should be expected if  solution $v$ to  original problem (\ref{mainequation})-(\ref{micv}) is smooth near the boundary $x_3=0$ near blowup time $T$. This scenario of the blowup can be called interior blowup. All other mild bounded ancient solutions then are related to the boundary blowup.  We note that the boundary blowup could still lead to a mild bounded ancient solution in the whole space. This might happen when the velocity tends to infinity fast in comparison with the rate at which $x$ is approaching the boundary.

 An interesting consequence of the above theorem is the following statement.
 \begin{pro}\label{smallness} There exists $\varepsilon>0$ (independent of $v$) such that if
 \begin{equation}\label{decayinx3}
|v(x,t)|\leq \frac \varepsilon{x_3}
\end{equation}
for all $x\in\mathbb R^3_+$ and $t\in ]0,T[$, the solution $v$ does not blow up.
\end{pro}

\pagebreak
\setcounter{equation}{0}
\section{Prelimaries}

Given $A<0$, consider the following Stokes problem in half space
$$\partial_tu-\triangle u+\nabla p=-{\rm div}\, F, \quad {\rm div}\,u=0$$
in $\mathbb R^3_+\times ]A,0[$,
$$u(x',0,t)=0$$
for $(x',t)\in\mathbb R^2\times ]A,0[$,
$$u(x,A)=u_0(x)$$
for $x\in\mathbb R^3_+$.

In addition, we assume that $u_0$ is a divergence free and $F$ and its first derivatives vanish at the boundary $x_3=0$. If $F$ and $u_0$ are sufficiently smooth and decay sufficiently fast at infinity, a solution of the above problem can be presented in the following way, see \cite{Sol1973},
\begin{equation}\label{maindecomposion}
    u=u^1+u^2,
\end{equation}
where
\begin{equation}\label{1stpart}
    u^1(x,t)=\int\limits_{\mathbb R^3_+}G(x,y,t-A)u_0(y)dy
\end{equation}
for any $x\in\mathbb R^3_+$ and for any $t\in ]A,0[$ and
\begin{equation}\label{2ndpartor}
    u^{2}(x,t)=-\int\limits^t_A\int\limits_{\mathbb R^3_+} G(x,y,t-\tau){\rm div}\,H(y,\tau)dy d\tau
\end{equation}
for the same $x$ and $t$,  where  $H=F+p^1\mathbb I$ and $p^1$ is a solution of the following Neumann boundary value problem
$$\Delta p^1=-{\rm div}{\rm div} F$$
in $R^3_+$ and $p^1_{,3}=0$ at $x_3=0$. Here, $G$  is
 the Green function for the Stokes system in the half space that has been studied in \cite{Sol1973} and has the form
\begin{equation}\label{GreenFunction}
    G=G^1+G^2,
\end{equation}
where
$$G^1_{ij}(x,y,t)=\delta_{ij}\Big(\Gamma(x-y,t)-\Gamma(x-y^*,t)\Big),$$
$$G^2_{i\beta}(x,y,t)=
4
\frac {\partial}{\partial x_\beta}\int\limits_0^{x_3}\int\limits_{\mathbb R^2}\frac{\partial E}{\partial x_i}(x-z)\Gamma(z-y^*,t)dz,\quad G^2_{i3}(x,y,t)=0.$$

As 
in the case of the whole space, it is convenient to present the function $u^2$
in the following equivalent way 

\begin{equation}\label{2ndpart}
    u^{2}_i(x,t)=\int\limits^t_A\int\limits_{\mathbb R^3_+}K_{ijm}(x,y,t-\tau)F_{jm}(y,\tau)dy d\tau
\end{equation}
for the same $x$ and $t$ as above.
The kernel $K$ has been introduced in \cite{Sol2003} and has the following structure
\begin{equation}\label{kernelK}
K_{ism}(x,z,t)=\overline{K}_{ism}(x,z,t)+\widehat{K}_{ism}(x,z,t),
\end{equation}
where $\overline{K}_{ism}(x,z,t)$ is a linear combination of the terms
$$\frac{\partial G_{ij}}{\partial z_k}(x,z,t)$$ and $\widehat{K}_{ism}(x,z,t)$ is a linear combination of the terms
$$\frac {\partial^2}{\partial x_\alpha\partial x_\beta}\int\limits_{\mathbb R^3_+}G_{ij}(x,y,t)\frac {\partial N^{(\pm)}}{\partial y_s}(y,z)dy.$$
Here, $N^{(\pm)}(x,y)=E(x-y)\pm E(x-y^*)$ with $y^*=(y',-y_3)$ and  $E(x)$ is fundamental solution to the Laplace equation in $\mathbb R^3$.

Let us outline how the above transformations can be done. Our arguments slightly differ from those that have been used in mentioned papers \cite{KochSol}
and \cite{Sol2003}. Consider the following boundary value problems
\begin{equation}\label{alpotential}
  \Delta_y \Phi_{mn}(x,y,t)=G_{mn}(x,y,t)
\end{equation}
with $\partial{\Phi_{mn}}/\partial{y_3}(x,y,t)=0$ if $n<3$ and with $\Phi_{mn}(x,y,t)=0$ if $n=3$ at $y_3=0$. And then integration by parts gives
$$K_{mjs}(x,y,t)=\frac{\partial^3 \Phi_{mj}}{\partial y_i\partial y_i\partial y_s}(x,y,t)-\frac{\partial^3 \Phi_{mn}}{\partial y_n\partial y_j\partial y_s}(x,y,t).$$
The above splitting for the potential $K$ can be easily derived with help of  the initial boundary value problems  for function $\Phi$, see (\ref{alpotential}), and some elementary properties of Green functions $G$, $N^+$, and $N^-$.

The following estimates for $G^i$  and $\widehat{K}$ have been obtained in papers \cite{Sol1973}, \cite{Sol2003}, and \cite{Sol2003UMN}:
$$ \Big|\frac {\partial^{|\alpha|+|\gamma|}G^2}{\partial x^\alpha\partial y^\gamma}(x,y,t-A)\Big|\leq c(\alpha,\gamma) (t-A)^{-\frac {\gamma_3}2}(t-A+x^2_3)^{-\frac{\alpha_3}2}\times $$\begin{equation}\label{GreenFunctEst2}
 \times  (|x-y^*|^2+t-A)^{-\frac {3+|\alpha'|+|\gamma'|}2}\exp{\Big(-\frac {cy^2_3}{t-A}\Big)},
\end{equation}
where $\alpha'=(\alpha_1,\alpha_2)$, $\gamma'=(\gamma_1,\gamma_2)$, and $|\gamma|=0,1$,
\begin{equation}\label{dergreen2}
    \Big|\frac{\partial G^1_{ij}}{\partial y_i}(x,y,t) \Big|+ |\widehat{K}_{ism}(x,y,t)|\leq \frac c{(|x-y|^2+t)^2}.
\end{equation}
\begin{equation}\label{derivativeintime}\Big|\partial^l_t G^2(x,y,t)\Big|\leq\frac c {t^l(|x'-y'|^2+x_3^2+y^2_3 +t)^\frac 32}\exp{\Big(-\frac {cy^2_3}{t}\Big)}. \end{equation}
for $l=0,1$.

Let  $K^1$ and $K^2$ be generated by $G^1$ and $G^2$, respectively.  In particular, we need the following estimate
\begin{equation}\label{hatK2}
  |\widehat{K}^2(x,y,t)|\leq \frac c{(|x-y^*|^2+t)^2}
\end{equation}
which can be obtained  
by an elementary modification of arguments used in the proof of Proposition 3.1 in \cite{Sol2003}.

\pagebreak
\setcounter{equation}{0}
\section{Scaling}

It is not difficult to show that there exists a sequence $(x^{(k)},t_k)$ with $x^{(k)}\in\mathbb R^3_+$ such that $t_k\to T-0$ and
\begin{equation}\label{blowupsequence}
    g(t_k)=M(t_k)=|v(x^{(k)},t_k)|\to\infty.
\end{equation}
\begin{remark}\label{remark on boundedness}From \cite{CKN} and from \cite{S3}, it follows that sequence $x^{(k)}$ is bounded. \end{remark}
We let $M_k=M(t_k)$. There are two main scenarios. In the first of the them,
\begin{equation}\label{scenario1}
    x^{(k)}_3M_k\to\infty
\end{equation}
and we scale $v$ and $q$ so that
\begin{equation}\label{Scaling}
    u^{(k)}(y,s)=\frac 1{M_k}v(x,t),\qquad p_k(y,s)=\frac 1{M^2_k}q(x,t),
\end{equation}
where
\begin{equation}\label{Sc1changeofvariables}
    y=M_k(x-x^{(k)}), \qquad s=M^2_k(t-t_k).
\end{equation}
By the above scaling, (\ref{mainequation}) and (\ref{mbcv}) are transformed into
\begin{equation}  \label{Scmainequation}  \partial_su^{(k)} +u^{(k)}\cdot\nabla u^{(k)}-\triangle u^{(k)}=-\nabla p_k ,\quad \mbox{div}\,u^{(k)}=0
\end{equation}
in $Q^k:= \{y=(y',y_3):\,\,y'\in\mathbb R^2,\,y_3>-x^{(k)}_3M_k\}\times  ]-t_kM^2_k,0 [$,
\begin{equation}\label{Sc1bcv}
    u^{(k)}(y',-x^{(k)}_3M_k,t)=0
\end{equation}
for any $y'\in\mathbb R^2$ and $s\in]-t_kM^2_k,0 [$. And, moreover, according to (\ref{blowupsequence}), we have
\begin{equation}\label{Sc1normirovka}
    |u^{(k)}(0)|=1.
\end{equation}

In the second scenario,
\begin{equation}\label{scenario2}
    x^{(k)}_3M_k\to a\in[0,\infty[.
\end{equation}
This suggests the same scaling  (\ref{Scaling}) but with slightly different change of variables
\begin{equation}\label{Sc2changeofvariables}
    y'=M_k(x'-x'^{(k)}),\quad y_3=M_kx_3, \qquad s=M^2_k(t-t_k),
\end{equation}
In this case, (\ref{mainequation}) and (\ref{mbcv}) are transformed into system (\ref{Scmainequation}), which is valid  in
 $R^3_+\times  ]-t_kM^2_k,0 [$, and into the boundary condition
\begin{equation}\label{Sc2bcv}
    u^{(k)}(y',0,t)=0
\end{equation}
for any $y'\in\mathbb R^2$ and $s\in]-t_kM^2_k,0 [$. Condition (\ref{Sc1normirovka}) is replaced with
\begin{equation}\label{Sc2normirovka}
    |u^{(k)}(0,x^{(k)}_3M_k,0)|=1.
\end{equation}

Our aim is to understand what happens if $k\to\infty$.

Without loss of generality, we may assume that the following statements are true:

\noindent
\textsc{Scenario 1} There exists a divergence free function $u\in L_\infty(Q_-)$  such that $|u|\leq 1$ a.e. in $Q_-$ and, for any $R>0$,
\begin{equation}\label{Scenario1conver}
  \int\limits_{Q(R)}u^{(k)}\cdot w dx dt\to \int\limits_{Q(R)}u\cdot w dx dt
\end{equation}
for any $w\in L_1(Q(R))$. Here, $Q(R)=B(R)\times ]-R^2,0[$;

\noindent
\textsc{Scenario 2}  There exists a divergence free function $u\in L_\infty(Q_-^+)$  such that $|u|\leq 1$ a.e. in $Q_-^+$ and, for any $R>0$,
\begin{equation}\label{Scenario2conver}
  \int\limits_{Q_+(R)}u^{(k)}\cdot w dx dt\to \int\limits_{Q_+(R)}u\cdot w dx dt
\end{equation}
for any $w\in L_1(Q_+(R))$. Here, $Q_+(R)=B_+(R)\times ]-R^2,0[$.

\pagebreak
\setcounter{equation}{0}
\section{Scenario 1}
Our goal is to show that the limit function $u$ from (\ref{Scenario1conver}) must be a mild bounded ancient solution  in the whole space  and satisfy the condition $$|u(0,0)|=1.$$

 Fix an arbitrary $A<0$. Let $k$ is sufficiently large so that $A>-t_kM^2_k$. We split $w:=u^{(k)}$ into two parts
$$w=w^1+w^2$$
so that
$$\partial_t w^1-\triangle w^1+\nabla r^1=0, \quad {\rm div}\,w^1=0$$
in $Q^k_A=\mathbb R^2\times\{y_3>-d_k\}\times ]A,0[$, where $d_k=x^{(k)}_3M_k\to\infty$,
$$w^1(y',-d_k,t)=0$$
for  $(y',t)\in\mathbb R^2\times ]A,0[$, and
$$w^1(y,A)=w(y,A)$$
for $y\in\mathbb R^2\times\{y_3>-d_k\}$.

The second part of $w$ is a solution to the following initial boundary value problem
$$\partial_t w^2-\triangle w^2+\nabla r^2={\rm div}\,F+\nabla p^{1(k)}, \quad {\rm div}\,w^2=0$$
in $Q^k_A$,
$$w^2(y',-d_k,t)=0$$
for  $(y',t)\in\mathbb R^2\times ]A,0[$, and
$$w^2(y,A)=0$$
for $y\in\mathbb R^2\times\{y_3>-d_k\}$.
Here, $F=w\otimes w$ and $p^{1(k)}$ is defined by the following Newmann boundary value problem
$$\triangle p^{1(k)}(y,t)=-{\rm div div}\, F(y,t), \qquad (y,t)\in Q^k_A,$$
$$p^{1(k)}_{,3}(y',-d_k,t)=0\qquad (y',t)\in \mathbb R^2\times ]A,0[.$$

Keeping in mind (\ref{GreenFunction}), we can present  $w^1$ in the  form
$$w^1=w^{1,1}+w^{1,2},$$
where
$$w^{1,1}(y,t)=\int\limits_{\mathbb R^2\times\{y_3>-d_k\}}G^1(y+d_k e_3,z+d_k e_3,t-A)w(z,A)dz$$
and
$$w^{1,2}(y,t)=\int\limits_{\mathbb R^2\times\{y_3>-d_k\}}G^2(y+d_k e_3,z+d_k e_3,t-A)w(z,A)dz.$$
Elementary calculations, estimate (\ref{GreenFunctEst2}), and the fact $|w(y,A)|\leq 1$ ensure an upper bound  for $w^{1,2}$
\begin{equation}\label{boundw12}
    |w^{1,2}(y,t)|\leq c\frac {\sqrt{t-A}}{y_3+d_k}.
\end{equation}

The similar decomposition can be exploited in order to evaluate  $w^2$. We have
$$w^2(y,t)=w^{2,1}+w^{2,2},$$
where
$$w^{2,1}(y,t)=$$$$=-\int\limits^0_A\int\limits_{\mathbb R^2\times\{y_3>-d_k\}}G^1(y+d_k e_3,z+d_k e_3,t-\tau)({\rm div}\, F+\nabla p^{1(k)})(z,\tau)dz d\tau$$
and
$$w^{2,2}(y,t)=$$$$=-\int\limits^0_A\int\limits_{\mathbb R^2\times\{y_3>-d_k\}}G^2(y+d_k e_3,z+d_k e_3,t-\tau)({\rm div}\, F+\nabla p^{1(k)})(z,\tau)dz d\tau$$
In fact, velocity  $w^{2,2}$ obeys the same estimate as (\ref{boundw12})
\begin{equation}\label{boundw22}
  |w^{2,2}(y,t)|\leq c\frac {\sqrt{t-A}}{y_3+d_k}.
\end{equation} To prove (\ref{boundw22}), we need the following lemma.
\begin{lemma}\label{w22} Let
$$V_i(x,t)=\int\limits^t_0\int\limits_{\mathbb R^3_+}K^2_{ijm}(x,y,t-\tau)H_{km}(y,\tau)dy d\tau$$
with $K^2_{ijm}=\overline{K}^2_{ijm}+\widehat{K}^2_{ijm}$,
where $\overline{K}^2_{ism}(x,z,t)$ is a linear combination of the terms
$$\frac{\partial G^2_{ij}}{\partial y_k}(x,z,t)$$ and $\widehat{K}^2_{ism}(x,z,t)$ is a linear combination of the terms
$$\frac {\partial^2}{\partial x_\alpha\partial x_\beta}\int\limits_{\mathbb R^3_+}G^2_{ij}(x,y,t)\frac {\partial N^{(\pm)}}{\partial y_s}(y,z)dy.$$
Then
$$|V(x,t)|\leq c\|H\|_\infty\frac t{x_3}.$$
\end{lemma}
\textsc{Proof} By (\ref{GreenFunctEst2}), we find 
$$|\overline{K}^2(x,y,t)|\leq \frac c{t^\frac 12(|x-y^*|^2+t)^\frac 32}\exp{\Big(-\frac {cy_3^2}t\Big)}.$$
For the second term $\widehat{K}^2$, we are going to make use of estimate (\ref{hatK2}). 
So, we have
$$|V(x,t)|\leq c\|H\|_{\infty}\int\limits_0^t\int\limits_{\mathbb R^3_+}|K^{2}(x,z,t-\tau)|dz d\tau\leq$$$$
\leq c\|H\|_{\infty}\int\limits_0^t\int\limits_{\mathbb R^3_+}\Big(\frac 1{\tau^\frac 12(|x-z^*|^2+\tau)^{\frac 32}}\exp{\Big(-\frac {cz_3^2}\tau\Big)}+\frac 1{(|x-z^*|^2+\tau)^2}\Big)dzd\tau\leq $$$$
 \leq c\|H\|_{\infty}\int\limits_0^t\int\limits_0^\infty\Big(\frac 1{\tau^\frac 12(|x_3+z_3|^2+\tau)^{\frac 12}}\exp{\Big(-\frac {cz_3^2}\tau\Big)}+\frac 1{|x_3+z_3|^2+\tau}\Big)dz_3d\tau\leq$$
$$\leq c\|H\|_{\infty}\Big(\frac 1{x_3}\int\limits_0^t\int\limits_0^\infty\frac 1{\tau^\frac 12}\exp{\Big(-\frac {cz_3^2}\tau\Big)}dz_3d\tau+\int\limits_0^t\int\limits_{x_3}^\infty\frac 1{|u|^2+\tau}du d\tau\Big)$$$$\leq c\|H\|_{\infty}\Big(\frac t{x_3}\int\limits_0^\infty\exp{(-cu^2)}du
+\int\limits_0^t\frac 1{\tau^\frac 12}\Big(\frac \pi 2-\arctan{\Big(\frac {x_3}{\tau^\frac 12}\Big)}\Big)d\tau\Big).
$$
To estimate the last integral, we do the following
$$\int\limits_0^t\frac 1{\tau^\frac 12}\Big(\frac \pi 2-\arctan{\Big(\frac {x_3}{\tau^\frac 12}\Big)}\Big)d\tau=\int\limits_0^t\frac 1{\tau^\frac 12}\arctan{\Big(\frac {\tau^\frac 12} {x_3}\Big)}d\tau \leq \frac t{x_3}.$$
Lemma \ref{w22} is proved. $\Box$

Now, upper bound (\ref{boundw22}) follows from Lemma \ref{w22} and the identity
$$w^{2,2}_i(y,t)=\int\limits^0_A\int\limits_{\mathbb R^2\times\{y_3>-d_k\}}K^2_{ijm}(y+d_k e_3,z+d_k e_3,t-\tau)F_{jm}(z,\tau)dz d\tau$$
with the potential $K^2$ derived from the boundary value problem (\ref{alpotential}), in which  $G$ is replaced with $G^2$.

If we let $w^0(y,t)=w^{1,1}+w^{2,1}$, then the new function
is a solution to the following initial boundary value problem
$$\partial_t w^0-\triangle w^0={\rm div}\,F+\nabla p^{1(k)},$$
in $Q^k_A$,
$$w^0(y',-d_k,t)=0$$
for  $(y',t)\in\mathbb R^2\times ]A,0[$, and
$$w^0(y,A)=w(y,A)$$
for $y\in\mathbb R^2\times\{y_3>-d_k\}$. Using exact representation formulae for $w^{1,1}$ and $w^{2,1}$, we may assume that $w^0$ is  bounded by a constant $c$, which is independent of $k$. On the other hand, we know that function $p^{1(k)}$, being extend to the whole $\mathbb R^3$ so that $p^{1(k)}(y',y_3,t)=p^{1(k)}(y',y_3+2d_k,t)$ for $y_3<-d_k$,
belongs to the space $L_\infty(A,0;BMO(\mathbb R^3))$ and the corresponding norm is bounded by a constant independent of $k$. So, we have
$$\sup\limits_{A\leq t\leq 0}\int\limits_{B(x,\sqrt{-A})}|p^{1(k)}(y,t)-[p^{1(k)}]_{B(x,\sqrt{-A})}(t)|^mdy\leq c(m)(-A)^\frac 32.$$
This means, see Appendix II, that sequence $w^0$ is precompact in  $C(K\times [A/2,0])$, where $K$ is an arbitrary compact  of $\mathbb R^3$. Now, it remains to make use of estimates for $w^{1,2}$ and $w^{2,2}$, pass to the limit in the equation for $w^0$ taking into account that $F=w\otimes w$ and conclude that, by arbitrariness of $A$, $u$ is a mild bounded ancient solution to the Navier-Stokes equations satisfying $|u(0,0)|=1$. $\Box$
\begin{remark}\label{convergScen1} In fact we have shown that $$u^{(k)}\to u$$
uniformly of on the closure of the set $Q(R)$ for any $R>0$.\end{remark}

\pagebreak
\setcounter{equation}{0}
\section{ Scenario 2}
Here, we are going to prove the following statement.
\begin{pro}\label{Sc2}
$$u^{(k)}\to u$$
uniformly of on the closure of the set $Q_+(R)=B_+(R)\times ]-R^2,0[$ with $B_+(R)=\{x\in B(R):\, x_3>0\}$ for any $R>0$. The limit function $u$ is  equal to zero at $y_3=0$ and is not trivial in the sense
$$|u(0,a,0)|=1$$
with $a>0$.
\end{pro}
\textsc{Proof} As in the  previous section, let us  split $w$ into two parts
$$w=w^1+w^2$$
where
$$\partial_t w^1-\triangle w^1+\nabla r^1=0, \quad {\rm div}\,w^1=0$$
in $Q_A=\mathbb R^3_+\times ]A,0[$,
$$w^1(y',0,t)=0$$
for  $(y',t)\in\mathbb R^2\times ]A,0[$, and
$$w^1(y,A)=w(y,A)$$
for $y\in\mathbb R^3_+$.

The second part is a solution to the following problem
$$\partial_t w^2-\triangle w^2+\nabla r^2={\rm div}\,F+\nabla p^{1(k)}, \quad {\rm div}\,w^2=0$$
in $Q_A$,
$$w^2(y',0,t)=0$$
for  $(y',t)\in\mathbb R^2\times ]A,0[$, and
$$w^2(y,A)=0$$
for $y\in\mathbb R^3_+$.
Here, $F=w\otimes w$ and $p^{1(k)}$ is defined by the following Neumann boundary value problem
$$\triangle p^{1(k)}(y,t)=-{\rm div div}\, F(y,t), \qquad (y,t)\in Q_A,$$
$$p^{1(k)}_{,3}(y',0,t)=0\qquad (y',t)\in \mathbb R^2\times ]A,0[.$$

Let us first discuss precompactness of $w^2$ in $C(\overline{B}_+(R)\times [A/2,0])$ for any positive $R$. Indeed, according to (\ref{GreenFunctEst2}), we have
$$    \Big|\frac{\partial G^2_{ij}}{\partial z_i}(x,z,t)\Big|\leq \frac c{t^\frac 12(|x-z^*|^2+t)^{\frac 32}}\exp{\Big(-\frac {cz_3^2}t\Big)}$$
and, using estimates (\ref{dergreen2}) and the definition of the kernel $K$, it is not difficult to show that
$$\int\limits_{\mathbb R^3_+}|K_{ism}(x,z,t)|dz d\tau\leq \frac c{\sqrt{t}}.$$

 Next, first, assuming $A\leq t_1<t_2\leq 0$ and $x^1,x^2\in \mathbb R^3_+$,
$$|w^{2}(x^1,t_1)-w^{2}(x^2,t_2)|\leq$$$$\leq\Big|\int\limits^{t_1}_A\int\limits_{\mathbb R^3_+}(K(x^1,y,t_1-\tau)-K(x^2,y,t_2-\tau))F(y,\tau)dy d\tau\Big|+$$$$+\Big|\int\limits^{t_2}_{t_1}\int\limits_{\mathbb R^3_+}K(x^2,y,t_2-\tau)F(y,\tau)dyd\tau\Big|=I_1+I_2.$$
For the second term, we have
$$I_2\leq c\|F\|_\infty\sqrt{t_2-t_1}.$$
Next, in the first term, we shall do the change of variables:
$$I_1\leq \|F\|_\infty\int\limits^{t_1-A}_0\int\limits_{\mathbb R^3_+}|K(x^1,y,\tau)-K(x^2,y,t_2-t_1+\tau)|dy d\tau\leq $$$$\leq \|F\|_\infty\int\limits^{-A}_0\int\limits_{\mathbb R^3_+}|K(x^1,y,\tau)-K(x^2,y,t_2-t_1+\tau)|dy d\tau=\|F\|_\infty J(x^1,x^2;t_1,t_2).$$
Using the above estimates of the Green function, it is not difficult to show that  given $\varepsilon>0$ there exists $\kappa(\varepsilon,A,R)$ such that
if  $|x^1-x^2|<\kappa$, $0<t_2-t_1<\kappa$, $x^1,x^2\in \overline{B}_+(R)$, and $A\leq t_1<t_2\leq 0$, then $J(x^1,x^2;t_1,t_2)<\varepsilon$.
So, the required precompactness of $w^2$ follows from Arzela-Ascoli theorem and from the following inequality
$$|w^{2}(x^1,t_1)-w^{2}(x^2,t_2)|\leq  \|F\|_\infty J(x^1,x^2;t_1,t_2)+ c \|F\|_\infty\sqrt{t_2-t_1}.$$

Precompactness of $w^1$ on the same sets is based on  similar arguments and the following fact
$$\int\limits_{\mathbb R^3_+}|G^2(x,y,t)|dy\leq c\sqrt{\frac t{x_3^2+t}}\leq c.$$

So, as usual, applying the diagonal Cantor procedure, we select a subsequence still denoted by $u^{(k)}$ that converges uniformly on sets $\overline{Q}_+(n)$ for any natural $n$. On the other hand,
$u^{(k)}(0,d_k,0)\to u(0,a,0)$ so that $|u(0,a,0)|=1$. The latter actually implies that $a>0$. $\Box$

The last arguments in the proof of Proposition \ref{Sc2} allow us to prove Proposition \ref{smallness}.

\textsc{Proof of Proposition \ref{smallness}}. If we assume that $T$ is a blow up time, then  since  that $a\leq A$ we have scenario 2 if one blows up our solution $v$ and gets a function $u$ of Proposition \ref{Sc2}.
It is not so difficult to see that modulus of continuity of $u$ in $\overline{Q}_+(2)$ depends only on  the integral modulus of continuity of the above Green's functions.  So,  there exists $\delta_0>0$ such that $|u(0,0,0)-u(0,a,0)|<1/2$ provided $a<\min(1,\delta_0)$ but in fact $|u(0,0,0)-u(0,a,0)|
=|u(0,a,0)|=1$. So, it remains to take $\varepsilon=\delta_0$. $\Box$

\begin{pro}\label{4p2} The limit function $u$ is a mild bounded ancient solution. \end{pro}

\textsc{Proof}  We have already proven that
\begin{equation}\label{unifconvergence}
    u^{(k)}\to u
\end{equation}
in $C(\overline{Q}_+(R))$.


We 
define $p^{1(k)}$ as a solution to the following Neumann boundary value problem
$$\triangle p^{1(k)}(x,t)=-{\rm div}\,{\rm div}\,\Big( u^{(k)}\otimes u^{(k)}\Big)(x,t)$$
for $(x,t) \in Q_-^+$ and
$$p^{1(k)}_{,3}(x',0,t)=0$$
for $(x',t)\in \mathbb R^2\times ]-\infty,0[$.
We extend $p^{1(k)}$ to the whole space $\mathbb R^3$ in the even way with respect to $x_3$ and $u^{(k)}$ is supposed to be extended by zero. Then the function $H^{(k)}:=u^{(k)}\otimes u^{(k)}+\mathbb I p^{1(k)}\in L_{\infty}(-\infty,0;BMO(\mathbb R^3))$ so that $$\sup\limits_k(\|H^{(k)}\|_{L_{\infty}(BMO)}+\|u^{(k)}\|_\infty)= d<\infty.$$ 


We know
$$u^{(k)}=\int\limits_{\mathbb R^3_+}G(x,y,t-A)u^{(k)}(y,A)dy+\int\limits^t_A\int\limits_{\mathbb R^3_+}K(x,y,t-\tau)F^{(k)}(y,\tau)dyd\tau,$$
where $F^{(k)}=u^{(k)}\otimes u^{(k)}$.
All the norms, bounded by $\|u^{(k)}\|_\infty$, $\|F^{(k)}\|_\infty$, or $\|H^{(k)}\|_{L_\infty(BMO)}$ only, remain to be bounded for the limit functions. 

We may transform  the above formula by integration by parts in the second term on the right hand side
$$u^{(k)}_i(x,t)=\int\limits_{\mathbb R^3_+}G_{ij}(x,y,t-A)u^{(k)}_j(y,A)d y+$$$$+\int\limits^t_A\int\limits_{\mathbb R^3_+}\frac {\partial G_{ij}}{\partial y_l}(x,y,t-\tau)(H^{(k)}_{jl}(y,\tau)-[H^{(k)}_{jl}]_{B((x',0),a)}(\tau))dy d\tau$$
where $B((x',0),a)$ is a ball of radius $a=(x_3^2+t-\tau)^\frac 12$ centered at the point $(x',0)$.

Then we split $u^{(k)}=u^1_{(k)}+u^2_{(k)}$ according to the decomposition of the Green function: $G=G^1+G^2$, see (\ref{GreenFunction}). And let $u^2_{(k)}=u^{2,1}_{(k)}+u^{2,2}_{(k)}$ so that
$$u^{2,2}_{i(k)}(x,t)=\int\limits^t_A\int\limits_{\mathbb R^3_+}\frac {\partial G^2_{ij}}{\partial y_l}(x,y,t-\tau)(H^{(k)}_{jl}(y,\tau)-[H^{(k)}_{jl}]_{B((x',0),a)}(\tau))dy d\tau.$$

From (\ref{GreenFunctEst2}), it follows
$$ \Big|\frac {\partial^{|\alpha|+|\gamma|}G^2}{\partial x^\alpha\partial y^\gamma}(x,y,t-A)\Big|\leq c(\alpha,\gamma) (t-A)^{-\frac {|\gamma|}2}(t-A+x^2_3)^{-\frac{|\alpha|}2}\times
$$\begin{equation}\label{GreenFunctionEstMod}\times  \exp{\Big(-\frac {cy^2_3}{t-A}\Big)}
{(|x'-y'|^2+y^2_3+x^2_3+t-A)^{-\frac {3}2}}.
\end{equation}
Now, we have
$$|\nabla^{|\alpha|}u^{2,1}_{(k)}(x,t)|\leq c (\alpha)(t-A+x^2_3)^{-\frac{|\alpha|}2}\times$$
$$\times \int\limits_{\mathbb R^3_+}(|x'-y'|^2+y^2_3+x^2_3+t-A)^{-\frac {3}2}\exp{\Big(-\frac {cy^2_3}{t-A}\Big)})|u^{(k)}(y,A)|dy$$
and
$$|\nabla^{|\alpha|}u^{2,2}_{(k)}(x,t)|\leq c (\alpha)\int\limits_A^t(t-\tau)^{-\frac 12}(t-\tau+x^2_3)^{-\frac{|\alpha|}2}\times$$
$$\times \int\limits_{\mathbb R^3_+}(|x'-y'|^2+y^2_3+x^2_3+t-\tau)^{-\frac {3}2}\exp{\Big(-\frac {cy^2_3}{t-\tau}\Big)})|H^{(k)}(y,\tau)-$$$$-[H^{(k)}]_{B((x',0),a)}(\tau)|dy d\tau.$$

 So,  we find
\begin{equation}\label{est21m}
  |\nabla^{|\alpha|}u^{2,1}_{(k)}(x,t)|\leq c (\alpha)d(x^2_3+t-A)^{-\frac {|\alpha|}2}\leq c (\alpha)d(x^2_3 -A)^{-\frac {|\alpha|}2}
\end{equation}
for all $(x,t)\in\mathbb R^3_+\times ]A/2,0[$ and, by Lemma \ref{integralestimate} with $\beta=\frac {t-\tau}c$ and $a^2=x^2_3+t-\tau$,
\begin{equation}\label{est22m}
  |\nabla^{|\alpha|}u^{2,2}_{(k)}(x,t)|\leq c (\alpha)d\int\limits_A^t(t-\tau)^{-\frac 12}(x^2_3+t-\tau)^{-\frac {|\alpha|}2}d\tau\leq $$
$$ \leq c (\alpha)d\int\limits^{\sqrt{t-A}}_0\frac {dq}{(x^2_3+q^2)^{\frac {|\alpha|}2}}\leq
c (\alpha)d\int\limits^{\sqrt{-A}}_0\frac {dq}{(x^2_3+q^2)^{\frac {|\alpha|}2}}
\end{equation}
for all $(x,t)\in\mathbb R^3_+\times ]A,0[$.

Next, 
we observe that $u^1_{(k)}$ is a solution to the following initial boundary value problem:
\begin{equation}\label{equatu1m}
  \partial_tu^1_{(k)}-\triangle u^1_{(k)}=-{\rm div}\,H^{(k)}
\end{equation}
for any $(x,t)\in\mathbb R^3_+\times ]A,0[$,
\begin{equation}\label{bcu1m}
  u^1_{(k)}(x',0,t)=0
\end{equation}
for any $(x',t)$ from $\mathbb R^2\times ]A,0[$,
\begin{equation}\label{inconu1h}
  u^1_{(k)}(x,A)=u^{(k)}(x,A)
\end{equation}
for all $x\in \mathbb R^3_+$.


Regarding $u^{2,2}_{(k)}(x,t)$, we can say the following: for an appropriate extension of it to the whole $\mathbb R^3$, we can claim that $\nabla u^{2,2}_{(k)}$ is $L_\infty(A,0;BMO)$ and
$\nabla u^{2,2}_{(k)}\in L_\infty(A,0;L_{q,{\rm unif}}(\mathbb R^3))$ for any finite $q\geq 1$.

Since $H^{(k)}$ is bounded in $L_\infty(A,0;BMO(\mathbb R^3))$, 
 by  estimates (\ref{est21m}) and (\ref{est22m}), we have (from the energy inequality)
\begin{equation}\label{fromenergy}
 \sup\limits_k \int\limits^0_{A/2}\int\limits_{\mathbb R^3_+}\varphi^2|\nabla u^1_{(k)}|^2dxdt<\infty
\end{equation}
for any test function $\varphi\in C^\infty_0(\mathbb R^3\times ]-A/2,A/2[)$.
Then, from the equation for the pressure $p^{1(k)}$, we derive the similar estimate for the gradient of $p^{1(k)}$
\begin{equation}\label{frompressure} \sup\limits_k\int\limits^0_{A/2}\int\limits_{\mathbb R^3_+}\varphi^2|\nabla p^{1(k)}|^2dxdt<\infty
\end{equation}
for any test function $\varphi\in C^\infty_0(\mathbb R^3\times ]-A/2,A/2[)$. Using coercive estimates for the heat equation, we then get
\begin{equation}\label{fromheat} \sup\limits_k\int\limits^0_{A/2}\int\limits_{\mathbb R^3_+}\varphi^2(|\partial_tu^1_{(k)}|^2+|\nabla^2 u^1_{(k)}|^2)dxdt<\infty
\end{equation}
for any test function $\varphi\in C^\infty_0(\mathbb R^3\times ]-A/2,A/2[)$.

Next, we may use a parabolic imbedding theorem to show that
 $$\sup\limits_k\int\limits^0_{A/2}\int\limits_{\mathbb R^3_+}\varphi^2|\nabla u^1_{(k)}|^ 3dxdt<\infty$$
 for the same test function $\varphi$. The same type of estimate is valid for the gradient of $p^{1(k)}$. Repeating these arguments several times, we can state that $\nabla u^1_{(k)}$ is bounded in domains $\mathbb R^3_+\times  ]3A/8,0[$.
  Boundedness of $\nabla p^{1(k)}$ follows from local regularity and the equation $\triangle p^{1(k)} =-u^{(k)}_{i,j}u^{(k)}_{j,i}$ being valid in $\mathbb R^3$ for appropriate extensions of $p^{1(k)}$ and $u^{(k)}_{i,j}u^{(k)}_{j,i}$ from $\mathbb R^3_+$ to $\mathbb R^3$. 
 Summarizing all these estimates, taking into account arbitrariness of $A$, and using a shift in time we find
 \begin{equation}\label{intermediateest2}
  \sup\limits_k |\nabla u^{(k)}(x,t)|\leq c\ln(2+x_3^{-1})
 \end{equation}
 for $(x,t)\in Q_-^+$ and 
\begin{equation}\label{intermediateest1}
  \sup\limits_k \sup\limits_{(x,t)\in Q^+_-}|\nabla p^{1(k)}(x,t)|\leq C<\infty.
\end{equation}

Now, we can go back to evaluation of function $u^2_{(k)}$. By (\ref{intermediateest2}) and (\ref{intermediateest1}), we do not need integrate by parts in the expression for $u^{2,2}_{(k)}$  any more and thus
$$u^{2,2}_{i(k)}(x,t)=-\int\limits^t_A\int\limits_{\mathbb R^3_+} G^2_{ij}(x,y,t-\tau)\Big[\frac {\partial}{\partial y_l}F^{(k)}_{jl}(y,\tau)+\frac {\partial}{\partial y_j}p^{1(k)}(y,\tau)\Big]dy d\tau.$$
 Next, according to  (\ref{intermediateest2}) and (\ref{intermediateest1}) and  by estimates  of Green function $G^2$, we find the following bound
 $$|\nabla u^{2,2}_{(k)}(x,t)|\leq c\int\limits^t_A\frac {d\tau}{(x_3^2+t-\tau)^\frac 12}\int\limits_0^\infty\frac {\exp{\Big(-\frac {cy^2_3}{t-\tau}\Big)}}{(x_3^2+y_3^2+t-\tau)^\frac 12}\ln(2+1/y_3)dy_3$$
 $$\leq c\int\limits^t_A\frac {d\tau}{(x_3^2+t-\tau)^\frac 34}\int\limits_0^\infty\frac {\exp{\Big(-\frac {cy^2_3}{t-\tau}\Big)}}{(x_3^2+y_3^2+t-\tau)^\frac 14}\ln(2+1/y_3)dy_3$$
 $$\leq  c\int\limits^t_A\frac {d\tau}{(x_3^2+t-\tau)^\frac 34}\Big[\int\limits^1_0 y_3^{-\frac 12}\ln(2+1/y_3)dy_3+\int\limits^\infty_1\exp{\Big(-\frac {cy^2_3}{t-\tau}\Big)}d y_3\Big]$$
 $$\leq c(1+\sqrt{-A})(-A)^\frac 14.$$
 The latter implies that
  in fact
\begin{equation}\label{intermediateest3}
 \sup\limits_k \sup\limits_{(x,t)\in Q^+_-}|\nabla u^{(k)}(x,t)|\leq C<\infty.
\end{equation}

Using estimates (\ref{GreenFunctEst2}), (\ref{GreenFunctionEstMod}), (\ref{intermediateest1}), and bootstrap arguments, we show
\begin{equation}\label{intermediateest}
  \sup\limits_{(x,t)\in Q^+_-, x_3\geq \delta}|\nabla^l u^{(k)}(x,t)|+|\nabla^{l+1} p^{1(k)}(x,t)|\leq C(l,\delta)<\infty
\end{equation}
for any $l\geq 0$ and for any $\delta>0$.

As to derivatives in $t$, by (\ref{derivativeintime}), we find 

$$|\partial_tu^{2,2}_{(k)}(x,t)|\leq c \int\limits^t_A\int\limits^\infty_0\int\limits_{\mathbb R^2}\frac 1{t-\tau}\frac 1{(|x'-y'|^2+x_3^2+y^2_3 +t-\tau)^\frac 32}\times $$$$\times \exp{\Big(-\frac {cy^2_3}{t-\tau}\Big)}\Big|{\rm div}\, H^{(k)}(y,\tau)\Big|dy d\tau\leq $$
$$\leq c \int\limits^t_A\frac 1{(t-\tau)^\frac 12}\frac 1{(x_3^2 +t-\tau)^\frac 12}d\tau.$$
So,
\begin{equation}\label{derintime22}
 \sup\limits_{x'\in\mathbb R^2} |\partial_tu^{2,2}_{(k)}(x',x_3,t)|\leq c \int\limits^{t-A}_0
 \frac 1{\vartheta^\frac 12}\frac 1{(x_3^2 +\vartheta)^\frac 12}d\vartheta
\end{equation}
 For the first term, we show in the same way that
$$|\partial_tu^{2,1}_{(k)}(x,t)|\leq c
\int\limits_{\mathbb R^3_+}|\partial_tG^2(x,y,t-A)||u^{(k)}(y,A)|dy $$ and thus \begin{equation}\label{derintime2}
|\partial_tu^{2,1}_{(k)}(x,t)|\leq \frac 1{(t-A)^\frac 12}\frac 1{(x_3^2 +t-A)^\frac 12}.\end{equation}
%
%

From equation (\ref{equatu1m}), 
(\ref{derintime22}), and (\ref{derintime2}), we can deduce that:
\begin{equation}\label{derintime}
  \sup\limits_{(x,t)\in Q^+_-, x_3\geq \delta}|\partial_tu^{(k)}(x,t)|\leq c(\delta).
\end{equation}
After passing to the limit as $k\to\infty$, we get that $u$ satisfies the standard integral identity with divergence free test functions. Then one can claim that there exists $p^2$ such that
\begin{equation}\label{pointwiseequation}
  \partial_tu-\triangle u+\nabla p^2=-{\rm div}\,H
\end{equation}
in $Q^+_-$.

Obviously, $p^2(\cdot,t) $ is a harmonic function in a half space whose gradient is bounded in $t$, in $x'$, and  in $x_3\geq \delta$ for any  $\delta>0$. Taking the limit in (\ref{equatu1m}), we show that $\nabla p^2$ satisfies
$$\partial_tu^2-\triangle u^2+\nabla p^2=0$$
in $\mathbb R^3_+\times ]A,0[$.  This, together with (\ref{est21m}), (\ref{est22m}) and  (\ref{derintime22}), (\ref{derintime2}), implies
$$\sup\limits_{x'\in\mathbb R^2}|\nabla p^2(x',x_3,t)|\to 0$$
as $x_3\to \infty$ for each $-\infty<t<0$. Moreover, more detail analysis of the above estimates for Green functions allows us to state:
$$|\nabla p^2(x,t)|\leq c\ln(2+
1/{x_3})$$
for all $(x,t)\in Q_-^+$. Proposition \ref{4p2} is proved. $\Box$


 \textsc{Proof of Theorem \ref{main}}  easily follows from the above arguments.

\pagebreak
\setcounter{equation}{0}
\section{Appendix I }

\begin{lemma}\label{integralestimate} Assume that numbers $m$ and $\alpha_0$ satisfy the condition
\begin{equation}\label{restriction}
  0<\alpha_0<\frac {m-1}3.
  \end{equation}
  Let $a$ and $\beta$ be positive numbers and let us define
  $$I(a,\beta):=\int\limits_{\mathbb R^3}\frac {|f(x)-[f]_{B(a)}|}{(|x|^2+a^2)^\frac 32}\exp{\Big(-\frac {x_3^2}\beta\Big)}dx,$$
where
  $[f]_{B(R)}$ is a mean value of $f$ over the ball $B(R)$ of radius $R$ centered at the origin. Let further $n=\frac m{m-1}$, $m_1=\frac m{1+\alpha_0}$, and $n_1=\frac {m_1}{m_1-1}$. Then
$$|I(a,\beta)|\leq c(m,\alpha_0)\|f\|_{BMO(\mathbb R^3)}a^{-\frac {3\alpha_0}m }\times $$
   \begin{equation}\label{intest}
   \times\Big(\sqrt {\beta}\int\limits^\infty_0\exp{(-u^2)}\frac {du}{(u^2\beta n_1/n+a^2)^{\frac {3n}{2n_1}-1}}\Big)^\frac 1n.
 \end{equation}
 \end{lemma}
  \textsc{Proof} Without loss of generality, we may assume that $[f]_{B(a)}=0$. We also let
  $$K(x):=\frac {1}{(|x|^2+a^2)^\frac 32}\exp{\Big(-\frac {x_3^2}\beta\Big)},\qquad A:=\|f\|_{BMO(\mathbb R^3)}.$$
  Then by H\"older inequality we have
  $$|I(a,\beta)|\leq\int\limits_{\mathbb R^3}K^{\frac 1{m_1}+\frac 1{n_1}}(x)|f(x)|dx\leq $$$$\leq \Big(\int\limits_{\mathbb R^3}K^\frac m{m_1}(x)|f(x)|^mdx\Big)^\frac 1m\Big(\int\limits_{\mathbb R^3}K^\frac n{n_1}(x)dx\Big)^\frac 1n.$$
For the second factor on the right hand side of the latter inequality, we find
$$\int\limits_{\mathbb R^3}K^\frac n{n_1}(x)dx=\int\limits^\infty_{-\infty}\exp{\Big(-\frac {x_3^2n}{\beta n_1}\Big)}
dx_32\pi\int\limits_0^\infty\frac {\varrho d\varrho}{(\varrho^2+x_3^2+a^2)^{\frac {3n}{2n_1}}}=$$
$$=2\pi\frac 1{\frac {3n}{2n_1}-1}\int\limits^\infty_{0}\exp{\Big(-\frac {x_3^2n}{\beta n_1}\Big)}
dx_3\frac 1{(x_3^2+a^2)^{\frac {3n}{2n_1}-1}}.$$
Here, we have used condition (\ref{restriction}) that ensures the inequality $\frac {3n}{2n_1}-1>0$. After changing variable, we have
$$\Big(\int\limits_{\mathbb R^3}K^\frac n{n_1}(x)dx\Big)^\frac 1n\leq c(m,\alpha) \Big(\sqrt {\beta}\int\limits^\infty_0\exp{(-u^2)}\frac {du}{(u^2\beta n_1/n+a^2)^{\frac {3n}{2n_1}-1}}\Big)^\frac 1n.$$

Next,
$$\int\limits_{\mathbb R^3}K^\frac m{m_1}(x)|f(x)|^mdx=\int\limits_{B(a)}K^\frac m{m_1}(x)|f(x)|^mdx+$$$$+\sum\limits^\infty_{k=0}\int\limits_{B(a2^{k+1})\setminus B(a2^k)}K^\frac m{m_1}(x)|f(x)|^mdx=I_1+I_2.$$
Obviously,
$$I_1\leq c(m)A^ma^{3-\frac {3m}{m_1}}.$$
For $I_2$, we are going to use the fact that $m/m_1=1+\alpha_0$ and two inequalities
$$|f|^m\leq 2^{m-1}(|f-[f]_{B(a2^{k+1})}|^m+|[f]_{B(a2^{k+1})}|^m)$$
and
$$|[f]_{B(a2^{k+1})}|\leq c(k+1)A.$$
The latter inequality can be found for example in \cite{Stein1993}, p. 141. Then, we have
$$I_2\leq c(m)\sum\limits^\infty_{k=0}\frac {|B(a2^{k+1})|}{(a2^k)^{3(1+\alpha_0)}}\frac 1{|B(a2^{k+1})|}\int\limits_{B(a2^{k+1})}|f-[f]_{B(a2^{k+1})}|^m+$$
$$+ c(m)\sum\limits^\infty_{k=0}\frac {|B(2^{k+1})|}{(a2^k)^{3(1+\alpha_0)}}|[f]_{B(a2^{k+1})}|^m\leq $$$$\leq c(m)A^m\sum\limits^\infty_{k=0}\frac 1{(a2^k)^{3\alpha_0}}(k+1)^m\leq c(m,\alpha_0)a^{-3\alpha_0}A^m.$$
This completes the proof of Lemma \ref{integralestimate}. $\Box$

\pagebreak
\setcounter{equation}{0}
\section{Appendix II}
Consider the following problem
$$\partial_tu-\Delta u=-{\rm div}\,g$$
in $\mathcal Q(a)=Q(a)+{(0,a^2)}$. Assuming that $u$ and $g$ are sufficiently smooth in $\mathcal Q(a)$, we wish to estimate the modulus of continuity of $u$ in the closure of the set $B(a/2)\times ]3(a/2)^2,a^2[$  in terms of $\|u\|_{\infty,\mathcal Q(a)}$ and $\|g-[g]_{B(a)}\|_{L_{m,\infty}(\mathcal Q(a))}$ for sufficiently large $m\geq 1$.

Fist we take a positive cut-off function $\varphi\in C^\infty_0(B(a)\times ]0,2a^2[)$
and extended it by zero. We let $w=\varphi u$. This function solves the following Cauchy problem
$$\partial_tw-\Delta w=-{\rm div}\,(\varphi f)+f\nabla\varphi+ w(\partial_t\varphi+\triangle\varphi)-2{\rm div}\,(u\nabla\varphi)$$
in $\mathbb R^3\times ]0,a^2[$, where $f(x,t)=g(x,t)-[g]_{B(a)}(t)$,
and $w(x,0)=0$ for $x\in \mathbb R^3$.

We split $w$ into two parts
$w=w^1+w^2$, where
$$\partial_tw^1-\Delta w^1=-{\rm div}\,(\varphi f)-2{\rm div}\,(u\nabla\varphi)$$
in $\mathbb R^3\times ]0,a^2[$ with $w^1(\cdot,0)=0$ in $ \mathbb R^3$ and
$$\partial_tw^2-\Delta w^2=f\nabla\varphi+ w(\partial_t\varphi+\triangle\varphi)$$
in $\mathbb R^3\times ]0,a^2[$ with $w^2(\cdot,0)=0$ in $ \mathbb R^3$. The second part can be estimated with the help of  coercive estimates
$$\|w^2\|_{W^{2,1}_m(\mathbb R^3\times ]0,a^2[)}\leq C(a,m,\|u\|_{\infty,\mathcal Q(a)},
\|f\|_{L_{m,\infty}(\mathcal Q(a))}).$$
For sufficiently large $m$,  they give an estimate for a H\"older norm of $u$ and thus we have a required estimate for the modulus of continuity of $w^2$.

For the first part, let us make use of the solution formula
$$w^1(x,t)=-\int\limits^t_0\int\limits_{\mathbb R^3}\Gamma(x-y,t-\tau)\Big({\rm div}\,(\varphi f)(y,\tau)+2{\rm div}\,(u\nabla\varphi))(y,\tau)\Big)dy d\tau$$
$$=-\int\limits^t_0\int\limits_{\mathbb R^3}\nabla_x\Gamma(x-y,t-\tau)\cdot\Big((\varphi f)(y,\tau)+2(u\nabla\varphi))(y,\tau)\Big)dy d\tau$$

We first assume that $m>2012$ and observe that the following simple estimate is true:
$$\int\limits_{\mathbb R^3}|K(x-y,t-\tau)|dy\leq \frac c{\sqrt{t-\tau}},$$
where $K(x-y,t)=\nabla_x\Gamma(x-y,t)$. Now, we wish to show that a given positive $\varepsilon$ there exists a positive $\delta(m,a)$ such that if $(x^1,t_1), (x^2,t_2)\in B(a)\times ]0,a^2[$ with $t_2>t_1$ and $|x^1-x^2|+t_2-t_1<\delta$, then
$$\int\limits_0^{a^2}\Big(\int\limits_{\mathbb R^3}|K(x^2-y,t_2-t_1+\vartheta)-K(x^1-y,\vartheta)|dy\Big)^\frac 1{m'}\frac 1{\vartheta^\frac 2m}d\vartheta<\varepsilon.$$
Assume that this is false. Then there exists $\varepsilon_0>0$ and sequences $$(x^{1,n},t_{1,n}), (x^{2,n},t_{2,n})\in B(a)\times ]0,a^2[$$ with $t_{2,n}>t_{1,n}$ and $|x^{1,n}-x^{2,n}|+t_{2,n}-t_{1,n}\to 0$ but
$$\int\limits_0^{a^2}\Big(\int\limits_{\mathbb R^3}|K(x^{2,n}-y,t_{2,n}-t_{1,n}+\vartheta)-K(x^{1,n}-y,\vartheta)|dy\Big)^\frac 1{m'}\frac 1{\vartheta^\frac 2m}d\vartheta\geq \varepsilon_0.$$
Indeed, for any positive $\vartheta$,
$$\int\limits_{\mathbb R^3}|K(x^{2,n}-y,t_{2,n}-t_{1,n}+\vartheta)-K(x^{1,n}-y,\vartheta)|dy\to 0$$
by the Lebesgue theorem.
And then
$$\int\limits_0^{a^2}\Big(\int\limits_{\mathbb R^3}|K(x^{2,n}-y,t_{2,n}-t_{1,n}+\vartheta)-K(x^{1,n}-y,\vartheta)|dy\Big)^\frac 1{m'}\frac 1{\vartheta^\frac 2m}d\vartheta\to 0,$$
since
$$\Big(\int\limits_{\mathbb R^3}|K(x^{2,n}-y,t_{2,n}-t_{1,n}+\vartheta)-K(x^{1,n}-y,\vartheta)|dy\Big)^\frac 1{m'}\frac 1{\vartheta^\frac 2m}\leq \frac c{\vartheta^\frac 1{2m'}}\frac 1{\vartheta^\frac 2m}\leq$$ $$\leq \frac c{\vartheta^{\frac 12 +\frac 3{2m}}}.$$
Since the function on the right hand side is integrable under our assumption on $m$, we arrive at a contradiction again by the Lebesgue theorem on the dominated convergence. The rest goes as follows:
$$|w^1(x^2,t_2)-w^1(x^1,t_1)|\leq $$$$\leq\Big|\int\limits^{t_2}_{t_1}\int\limits_{\mathbb R^3}K(x_2-y,t_2-\tau)\cdot\Big((\varphi f)(y,\tau)+2(u\nabla\varphi))(y,\tau)\Big)dy d\tau\Big|+$$
$$+\Big|\int\limits^{t_1}_{0}\int\limits_{\mathbb R^3}(K(x_2-y,t_2-\tau)-K(x_1,t_1-\tau))\cdot\Big((\varphi f)(y,\tau)+2(u\nabla\varphi))(y,\tau)\Big)dy d\tau\Big|.$$
Then we apply H\"older inequality
and have
$$|w^1(x^2,t_2)-w^1(x^1,t_1)|\leq $$$$\leq\int\limits^{t_2}_{t_1}\Big(\int\limits_{\mathbb R^3}\Big|K(x_2-y,t_2-\tau)\Big|dy\Big)^\frac 1{m'}\times$$$$\times\Big(\int\limits_{B(a)}\Big|K(x_2-y,t_2-\tau)\Big|(|f|^m+|u|^m)dyd\tau\Big)^\frac 1m+$$
$$+\int\limits^{a^2}_{0}\Big(\int\limits_{\mathbb R^3}\Big|K(x_2-y,t_2-\tau)-K(x_1,t_1-\tau)\Big|dy\Big)^\frac 1{m'}\times$$$$\times\Big(\int\limits_{B(a)}\Big|K(x_2-y,t_2-\tau)-
K(x_1,t_1-\tau)\Big|(|f|^m+|u|^m)dyd\tau\Big)^\frac 1m\leq $$
$$\leq C(a,m,\|u\|_{\infty,\mathcal Q(a)},\|f\|_{L_{m,\infty}(\mathcal Q(a))})\Big(\int\limits^{t_2}_{t_1}\frac 1{(t_2-\tau)^{\frac 12 +\frac 3{2m}}}d\tau+$$
$$+\int\limits_0^{a^2}\Big(\int\limits_{\mathbb R^3}|K(x^2-y,t_2-t_1+\vartheta)-K(x^1-y,\vartheta)|dy\Big)^\frac 1{m'}\frac 1{\vartheta^\frac 2m}d\vartheta\Big)=$$
$$=C(a,m,\|u\|_{\infty,\mathcal Q(a)},\|f\|_{L_{m,\infty}(\mathcal Q(a))})\Big((t_2-t_1)^{\frac 12 -\frac 3{2m}}d+$$
$$+\int\limits_0^{a^2}\Big(\int\limits_{\mathbb R^3}|K(x^2-y,t_2-t_1+\vartheta)-K(x^1-y,\vartheta)|dy\Big)^\frac 1{m'}\frac 1{\vartheta^\frac 2m}d\vartheta\Big).$$
This, together with the above statements, gives a required estimate for the modulus of continuity in the closure of the $B(a/2)\times ]3(a/2)^2,a^2[$.

\pagebreak

\end{document}